\documentclass[a4paper,12pt]{article}
\usepackage{amsthm}
\usepackage{amssymb}
\usepackage{amsmath}
\usepackage{amsfonts}

\newtheorem{remark}{Remark}[section]

\newtheorem{lemma}{Lemma}[section]
\newtheorem{theorem}{Theorem}[section]

\def\b1{\mbox{\boldmath $1$}}

\newenvironment{demo*}{\vspace{3mm}\noindent{\bf Proof.}}{\hfill $\Box$ \vspace{3mm}}

\begin{document}

\baselineskip=20pt
\title{\bf \normalsize CONVEXITY OF RUIN PROBABILITY  AND OPTIMAL DIVIDEND \\
 STRATEGIES FOR A GENERAL L\'EVY PROCESS}
\author{\normalsize$^a${\sc Chuancun Yin},  $^b${\sc Kam Chuen Yuen} and  $^a${\sc Ying Shen} \\
{\normalsize\it $^a$School of Mathematical Sciences, Qufu Normal University}\\
{\normalsize\it Shandong 273165, P.R.\ China} \\
 e-mail: ccyin@mail.qfnu.edu.cn\\
[3mm] {\normalsize\it $^b$Department of Statistics and Actuarial
Science, The University of Hong Kong}\\
\noindent{\normalsize\it Pokfulam Road, Hong Kong}\\
e-mail: kcyuen@hku.hk } \maketitle
 \centerline{\large {\bf Abstract}}
In this paper, we  consider the   optimal dividends  problem for a company whose cash reserves follow  a general L\'evy process  with  certain positive jumps and arbitrary
negative jumps. The objective is to find a policy which maximizes the expected  discounted dividends until the  time of ruin. Under appropriate conditions,  we appeal to very recent results in the theory of potential analysis of subordinators to obtain the convexity  properties of probability of ruin. We present conditions under which the optimal dividend strategy,   among all admissible ones, takes the form   of a barrier strategy.
\vspace{1cm}



\normalsize
\baselineskip=18pt
\section{\normalsize Introduction}\label{intro}
In the literatures of actuarial science and finance, the optimal dividend problem is one of the key topics. For companies paying dividends to shareholders, a commonly-encountered problem is to find a dividend strategy that maximizes the expected total discounted dividends until ruin.  The pioneer work can be traced to  de Finetti (1957) who considered a discrete-time risk model with step sizes $\pm 1$ and showed that a certain barrier strategy maximizes the expected discounted dividend payments. Since then, the problem of finding the optimal dividend strategy has become a popular topic in the actuarial literature. For diffusion models, see, for example,   Jeanblanc-Picqu\'e and Shiryaev (1995),  Asmussen and Taksar (1997),  Paulsen (2003),    Gerber and Shiu (2004), L{\o}kka and Zervos (2008), Paulsen (2008), He and Liang (2009). Bai and Paulsen (2012).
For the Cram\'er-Lundberg risk model, some related works on this subject include, among others, Gerber (1969), Azcue and Muler (2005, 2010), Yuen et al. (2007),   Kulenko and Schimidli (2008), Bai and Guo (2010) and Hunting and   Paulsen (2013).

Analysis of optimal dividends for L\'evy risk processes are of particular interest which have undergone an intensive development.  For example, Avram et al. (2007) considered a general spectrally negative L\'evy process and gave a sufficient condition involving the generator of the L\'evy process for the optimality of barrier strategy; Loeffen (2008) showed that barrier strategy is optimal among all admissible strategies for general spectrally negative L\'evy risk process with completely monotone jump density; Kyprianou, Rivero and Song (2010) relaxed this condition on the jump density; Yin and Wang (2009) also studied the same problem and gave an alternate proof of the result; Loeffen (2009a, 2009b) considered the optimal dividend problem with transaction costs and a terminal value for the spectrally negative L\'evy process. Recently, Bayraktar, Kyprianou and Yamazaki (2013) using the fluctuation theory of  spectrally positive  L\'evy processes, show  the optimality of barrier strategies for all such L\'evy processes. See Yin and Wen (2013) for a different approach. All of the above mentioned works are based on spectrally one-sided models.   There are, however,  few papers studied the analogous problems for  L\'evy process with two-sided jumps (cf. Bo et al. (2012)). Inspired by the works of Avram et al. (2007), Loeffen (2008) and Kyprianou et al. (2010), Yuen and Yin (2011) considered the optimal dividend problem for a special L\'evy process with both upward and downward jumps and showed that the optimal strategy takes the form of a barrier strategy if the L\'evy measure (both negative and positive jumps) has a completely monotone density.  The purpose of the present paper is to extend the result of Yuen and Yin (2011) to the case with less restrictive conditions on the L\'evy measure.  Although the broader case definitely makes the optimization problem more challenging and complex, recent results on the theory of potential analysis of subordinators can be applied to handle it. Specifically, our main results show that the optimal dividend strategy is still of a barrier type if the L\'evy process  has  certain positive jumps  and L\'evy density of negative jumps is completely monotone or log-convex.

The paper is organized as follows. In Section 2, we introduce the mathematical formulation of the problem. In Section 3, we give a brief review on ladder processes and potential measure for general L\'evy processes.  The convexities of the ruin probability and the  scale function are discussed in Sections 4 and 5 and the main results and their proofs are given in Section 6.
\section{\normalsize The model}\label{math}
To present the mathematical formulation of the problem of study, let us first introduce some notations and definitions.
Let $X = \{X_t\}_{t\ge 0}$ be a real-valued  L\'evy process on a filtered probability space $(\Omega, \cal{F}, \Bbb{F}, P)$ where $\Bbb{F}=({\cal F}_t)$$_{t\ge 0}$ is generated by the process $X$ and  satisfies the usual conditions of right-continuity
and completeness.  Denote by $P_x$ for the  law of $X$  when $X_0=x$. Let $E_x$ be the expectation  associated with   $P_x$. For notational convenience, we write $P$ and $E$ when $X_0=0$.
Write the L\'evy triplet of $X$ as $(a, \sigma^2, \Pi)$, where $a,\sigma\ge 0$ are real constants and $\Pi$ is a positive measure on $(-\infty,\infty)\setminus \{0\}$ which satisfies the integrability condition
$$\int_{-\infty}^{\infty}(1 \wedge x^2)\Pi(dx)<\infty.$$
  If $\Pi(dx) =\pi(x)dx$, then we call $\pi$ the L\'evy density. The characteristic exponent of $X$ is given by
$$\kappa(\theta)=-\frac{1}{t}\log E(e^{i\theta X_t}) =-ia\theta+\frac1 2\sigma^2\theta^2 +\int_{-\infty}^{\infty}(1-e^{i\theta x}+i\theta x\text{\bf
1}_{\{|x|<1\}})\Pi (dx),$$
where $\text{\bf 1}_{A}$ is the indicator of set $A$. Furthermore, define the Laplace exponent of $X$ by
\begin{equation}
  \Psi (\theta)=\frac1{t}\log E(e^{\theta X_t}) =a\theta + \frac1 2\sigma^2\theta^2
+\int_{-\infty}^{\infty}(e^{\theta x}-1-\theta x\text{\bf
1}_{\{|x|<1\}})\Pi (dx). \label{math-eq1}
\end{equation}
Such a L\'evy process is of bounded variation if and only if $\sigma=0$ and $\int_{-1}^{1}|x|\Pi(dx)<\infty$. If $\Pi\{(0,\infty)\}=0$, then the L\'evy process $X$ with no positive jumps is called the spectrally negative L\'evy process;   If $\Pi\{(-\infty,0)\}=0$, then the L\'evy process $X$ with no negative jumps is called the spectrally positive L\'evy process. It is usual to assume  that $P(\lim_{t\to\infty}X_t=+\infty)=1$ which says nothing other than $\Psi'(0+)>0$. For more information on L\'evy processes we refer to  the excellent book by Kyprianou (2006).

Now, we consider an insurance company or investment company whose cash reserve process (also called risk process or surplus process) evolves according
to the process $X$ before dividends are deducted. Let $ \xi=\{L_t^{ \xi}:t\ge 0\}$ be a dividend policy  consisting of
a right-continuous non-negative non-decreasing process adapted to the filtration $\{\cal{F}$$_t$\}$_{t\ge 0}$ of $X$ with $L_{0-}^{ \xi}=0$,
where $L_t^{ \xi}$ represents the cumulative dividends paid up to time $t$. Given a control policy $\xi$, the controlled reserve process with initial capital $x\ge 0$ is given by $U^{\xi}=\{U_t^{\xi}:t\ge 0\}$ where
\begin{equation}
U_t^{\xi}=X_t-L_t^{ \xi},\label{math-eq2}
\end{equation}
with $X_0=x$. Let $\tau^{\xi}=\{t>0: U_t^{ \xi}<0\}$ be the ruin time when dividend payments are taken into account. Define the value function associated to dividend policy $\xi$ by
$$V_{\xi}(x)=E_x\left(\int_0^{\tau^{ \xi}}\text{e}^{-\delta
t}dL_t^{ \xi}\right),$$
where $\delta>0$ is the discounted rate. The integral is understood pathwise in a Lebesgue-Stieltjes sense. Clearly, $V_{\xi}(x)=0$ for $x<0$. A dividend policy is called admissible if $L_{t}^{ \xi}-L_{t-}^{ \xi}\le U_t^{ \xi}$ for $t<\tau^{
\xi}$ and $L_{\tau^{\xi}}^{ \xi}-L_{\tau^{\xi}-}^{\xi}=0$ for $\tau^{\xi}<\infty$. Denote by $\Xi$ the set of all admissible dividend policies. Our objective is to find $$V_*(x)=\sup_{ \xi\in\Xi}V_{ \xi}(x),$$
and an optimal policy $ \xi^*\in\Xi$ such that $V_{\xi^*}(x)=V_*(x)$ for all $x\ge 0$. The function $V_*$ is called the optimal value function.

We denote by $\xi_b=\{L_t^b:t\ge 0\}$ the barrier strategy at $b$ and let $U^b$ be the corresponding risk process, that is, $U_t^b=X_t-L_t^b$. Note that $\xi_b\in \Xi$. Also, if $U^b_0\in [0,b]$, then the process $L^b_t$ can be explicitly represented by
$$L^b_t=(\sup_{s\le t}X_s-b)\vee 0.$$
 If $U_0^b=x>b$, then
$$L^b_t=(x-b)\text{\bf 1}_{\{t=0\}}+(\sup_{s\le t}X_s-b)\vee 0.$$
Denote by $V_b(x)$ the dividend-value function if barrier strategy $\xi_b$ is applied, that is,
\begin{equation}
V_b(x)=E_x\left(\int_0^{\tau^{\xi_b}}e^{-\delta t}dL_t^b\right).  \label{math-eq3}
\end{equation}
Applying Ito's formula for semimartingale, we can prove that $V_b$ is the solution to
 $$
  \left\{
  \begin{array}{ll}
  \Gamma
V_b(x)=\delta V_b(x),\; & x>0,\\
     V_b(x)=0, &x<0,\\
      V_b(0)=0, &\sigma^2>0,\\
        V_b'(b)=1,\\
     V_b(x)=x-b+V_b(b),
  \end{array}
  \right.
  $$
where $\Gamma$ is the infinitesimal generator of $X$ with
\begin{equation}
\Gamma g(x)=\frac{1} {2}\sigma^2 g''(x)+a g'(x)
+\int_{-\infty}^{\infty}[g(x+y)-g(x)-g'(x)y\text{\bf
1}_{\{|y|<1\}}]\Pi(dy).\label{math-eq4}
\end{equation}
In the sequel, we assume that,    for any $\delta>0$,  the
equation $\Psi(z) =\delta$ has a unique  solution on $(0,\infty)$, say  $\rho(\delta)$.
A typical example  is that the L\'evy measure of the
 positive jumps has the following  gamma
distribution $\Gamma(r,1/\gamma)$, i.e.,
$$P(x)=\int_0^x\frac{r^{1/\gamma}}{\Gamma(1/\gamma)}y^{1/\gamma-1}\text{e}^{-r y}
\text{d}y, \;x>0,$$
 where   $r$ is a positive number and $\gamma$ is an even  number.

Following similar reasoning to  Yuen and Yin (2011), $V_b$ can be expressed as
\begin{equation}
  V_b(x)=\left\{
  \begin{array}{ll}
    \frac{h(x)}{h'(b)},&0\le x\le b,\\
    x-b+\frac{h(b)}{h'(b)}, &x>b,
  \end{array}
  \right.\label{math-eq5}
\end{equation}
where
\begin{equation}
h(x)=[1-\tilde{\psi}(x)]e^{\rho(\delta) x}.\label{math-eq6}
\end{equation}
Here, $\tilde{\psi}(u)$ be the ruin probability for a L\'evy process $\tilde{X}$ with Laplace exponent
$\psi_{\rho(\delta)}$ given by $\psi_{\rho(\delta)}(\eta)=\Psi(\eta+\rho(\delta))-\delta$.
Note that the process $\tilde{X}$ has the L\'evy triplet
$(\tilde{a}, \tilde{\sigma}^2, \tilde{\Pi})$,
 where $\tilde{\sigma}^2=\sigma^2$, $\tilde{\Pi}(dx)=e^{\rho(\delta) x}\Pi(dx),$ and
$$\tilde{a}=a+\sigma^2\rho(\delta)+\int_{-\infty}^{k}
(e^{\rho(\delta) y}-1)y\text{\bf 1}_{\{|y|\le 1\}}\Pi(dy).$$
Moreover, \begin{equation} \int_{|x|\ge 1}e^{\rho(\delta)
x}\Pi(dx)<\infty.  \nonumber
\end{equation}


 \vskip 0.2cm
\section{\normalsize Some results on ladder processes and potential measure}\label{ladd}
\setcounter{equation}{0}

In this section, we recap some basic facts about ladder processes and potential measure.
Consider  the dual process $Y=\{Y_t\}_{t\ge 0}$, with $Y_0=0$, where
$Y_t=-X_t$, $t\ge 0$. It is easy to see that that the  L\'evy triplet
of $Y$ is  $(-a, \sigma^2, \Pi_Y)$, where $\Pi_Y(dx)=\pi_X(-x)dx$.
 Let
$$\underline{Y_{t}}=\inf_{0\le s\le t}Y_s\;\;\text{and}\; \overline{Y}_{t}=\sup_{0\le s\le t}Y_s,$$
be  the processes of the first infimum and the last supremum
 of the L\'evy process $Y$, respectively. Following
Kl\"uppelberg, Kyprianou and Maller (2004), we now introduce the notion
of ladder processes and potential measure.
 Let $L =\{L_t :t\ge
0\}$ denote the local time in the time period $[0, t]$ that
$\overline{Y}-Y$ spends at zero. Then $L^{-1}=\{L^{-1}_t: t\ge 0\}$
is the inverse local time such that $L_t^{-1}=\inf\{s\ge 0: L_s>t\}$,
where we take the infimum of the empty set as $\infty$. Define an
increasing process $H$  by $\{H_t =Y_{L_t^{-1}}: t\ge 0\}$, that is, the
process of new maxima indexed by local time at the maximum. The
processes $L^{-1}$ and $H$ are both  defective subordinators, and we
call them the ascending ladder time and ladder height process of
$Y$, respectively. It is understood that $H_t=\infty$ when $L^{-1}_t=\infty$.
Throughout the paper, we denote the nondefective versions of
$L, L^{-1}, H$ by $\cal{L}$, ${\cal{L}}^{-1}$, $\cal{H}$ ,
respectively. In fact, the pair (${\cal{L}}^{-1}$, $\cal{H})$ is a
bivariate subordinator. Define $(\hat{L}^{-1}, \hat{H})$
the descending ladder time and the ladder height processes in an
analogous way. Note that $\hat{H}$ is a process which is
negatively valued. Because $Y$ drifts to $-\infty$,  the
decreasing ladder height process is not defective.  Associated
with the ascending and descending ladder processes are the bivariate
renewal functions $U$ and $\hat{U}$. The former is defined by
$$U(dx,ds) =\int_0^{\infty}P(H_t\in dx,L^{-1}_t\in ds)dt.$$
Taking Laplace transforms shows that
$$\int_0^{\infty}\int_0^{\infty}e^{-\beta x-\alpha
s}U(dx,ds)=\frac{1}{k(\alpha,\beta)}, \ \ \quad \ \alpha,\beta\ge
0,$$
where $ k(\alpha,\beta)$ is its joint Laplace exponent such that
$$k(0, \beta) = q + c\beta + \int_{(0,\infty)}(1-e^{-\beta x})\Pi_H(dx),$$
$q\ge 0$ is the killing rate of $H$ so that $q >0$ if and only if $\lim_{t\to\infty}Y_t =-\infty$, $c
\ge 0$ is the drift of $H$, and $\Pi_H$ is its jump measure. Denote the marginal measure of $U(\cdot,\cdot)$ by
\begin{equation}
U(dx)=U(dx,[0,\infty))= \int_0^{\infty}P(H_t\in
dx)dt=\int_0^{\infty}e^{-q t}P({\cal H}_t\in dx)dt,\; x\ge
0.\label{1add-eq1}
\end{equation}
The function $U$ is called the potential/renewal measure. As for the descending ladder process, $\hat{U}$ and $\hat{k}$ are defined similarly. Write $\Pi_+$ and $\Pi_-$ for the restrictions of $\Pi(du)$ and
$\Pi(-du)$ to $(0,\infty)$. Furthermore, for $u>0$, define
$$\bar{\Pi}^{+}_{Y}(u)=\Pi_Y\{(u,\infty)\},\ \ \bar{\Pi}^{-}_Y(u)=\Pi_Y\{(-\infty,-u)\},
\ \ \bar{\Pi}_Y(u)=\bar{\Pi}_Y^+(u)+\bar{\Pi}_Y^-(u).$$

We next introduce the notions of a special Bernstein function and complete Bernstein function and two useful results.
Recall that a function $\phi:(0,\infty)\rightarrow (0,\infty)$ is
called a  Bernstein function if it admits a representation
$$\phi(\lambda)=a+b\lambda+\int_0^{\infty}(1-e^{-\lambda x})\mu(dx),$$
where $a\ge0$ is the killing term, $b\ge 0$ is the drift, and $\mu$ is the L\'evy measure concentrated on $(0,\infty)$ satisfying
$\int_0^{\infty}(1\wedge x)\mu(dx)<\infty$.  A function $\psi$ is called a special Bernstein function if the
function $\psi(\lambda)= \lambda/\phi(\lambda)$ is again a
Bernstein function. Let
$$\psi(\lambda)=\tilde{a}+\tilde{b}\lambda+\int_0^{\infty}(1-e^{-\lambda x})\nu(dx)$$ be the corresponding representation. It was shown in Song and Vondrac$\check{e}$k (2006) that
$$\tilde{b}=\frac{1}{a+\mu((0,\infty))}\text{\bf 1}_{\{b=0\}},\;\; \tilde{a}=\frac{1}{b+\int_0^{\infty}t\mu(dt)}\text{\bf 1}_{\{a=0, \mu(0,\infty)<\infty\}}.$$
A possibly killed subordinator is called a special subordinator if
its Laplace exponent is a special Bernstein function. Song and Vondra$\check{c}$ek (2010) showed that a sufficient condition
for $\phi$ to be a special subordinator is that $\mu(x,\infty)$ is
log-convex on $(0,\infty)$.  A function $\phi: (0,\infty)\rightarrow \Bbb{R}$ is called a
complete Bernstein function if there exists a Bernstein function
$\eta$ such that
$$\phi(\lambda)=\lambda^2\cal{L}\eta(\lambda), \ \ \lambda> \text{0},$$
where $\cal{L}$ stands for the Laplace transform. It is known
that every complete
Bernstein function is a Bernstein function and that the following
three conditions are equivalent:
\begin{itemize}
\item[(i)] $\phi$ is a complete Bernstein function;
\item[(ii)] $\psi(\lambda)= \lambda/\phi(\lambda)$ is a complete Bernstein function;
\item[(iii)] $\phi$ is a Bernstein function whose L\'evy measure $\mu$ is
given by $$\mu(dt)=dt\int_0^{\infty}e^{-st}\nu(ds),$$
\end{itemize}
where $\nu$ is
a measure on $(0,\infty)$ satisfying
$$\int_0^1\frac{1}{s}\nu(ds)<\infty,\ \
\int_1^{\infty}\frac{1}{s^2}\nu(ds)<\infty.$$
To end the section, we present two results which are useful in potential theory and will be used in later sections
of the paper. The first due to Kyprianou, Rivero and Song (2010) (see also Song and  Vondracek (2010)) is summarized in Lemma 3.1 while the second due to Kingman (1967) and Hawkes (1977) is given in Lemma 3.2.
\begin{lemma}\label{ladd-1}
Let $H$ be a subordinator whose L\'evy density, say ${\mu}(x), x>0$,
is log-convex.  Then, the restriction of its potential measure to
$(0,\infty)$   has a non-increasing and convex density. Furthermore, if the drift of $H$ is strictly positive, then the
density is in  $C^1(0,\infty)$.
\end{lemma}

\begin{lemma} \label{ladd-2}
Suppose that $H$ is a subordinator with Laplace exponent $\phi$ and potential measure $U$. Then, $U$ has a density $u$ which is
completely monotone on $(0,\infty)$ if and only if the tail of the
L\'evy measure is  completely monotone.
\end{lemma}

\begin{remark}
Note that the tail of the L\'evy measure $\mu$ is a completely
monotone function if and only if $\mu$ has a completely monotone
density. Thus, we have the following two equivalent statements: $\phi$ is a
complete Bernstein function if and only if $U$ has a density $u$
which is completely monotone on $(0,\infty)$; or equivalently, $U$
has a density $u$ which is completely monotone on $(0,\infty)$ if
and only if $\mu$ has a completely monotone density.
\end{remark}
\setcounter{equation}{0}
\section{\normalsize Convexity of probability of ruin }\label{ruin}
\setcounter{equation}{0}
Define the probability of ruin by
\begin{eqnarray*}
\psi(x)&=&P({\rm there\ exists} \  t\ge 0\ {\rm such\ that}\ x+X_t\le
0)\\
&=&P({\rm there\   exists}  \ t\ge 0  \ {\rm  such \ that}  \ Y_t\ge
x).\end{eqnarray*} It follows from Bertion and Doney (1994) that
$\psi(x)=\alpha U(x,\infty)$, where
$\alpha^{-1}=U{(0,\infty)}=\int_0^{\infty}P(H_t<\infty)dt$, with $U$
given in (\ref{1add-eq1}).

For simplicity, we write the L\'evy measure $\Pi$ as
$$\Pi(dx)=\left\{\begin{array}{ll} \Pi_{+}(dx), &  x>0,\\
 \pi_{-}(-x)dx, & x<0. \end{array}\right.$$
Recall that an infinitely differentiable function $f\in
(0,\infty)\rightarrow [0,\infty)$  is called completely monotone if
$(-1)^n f^{(n)}(x)\ge 0$ for all $n=0,1,2,\cdots$ and all $x>0$.
\begin{lemma}\label{ruin4-1} \ (Vigon (2002)) \ For the L\'evy process $X$, we have
$$\overline{\Pi}_{\cal{H}}(x) =-\int_{-\infty}^0\hat{U}(dy) {\overline \Pi}_Y^+(x-y)
=-\int_{-\infty}^0\hat{U}(dy) {\overline \Pi}_X^{\ -}(x-y),\; \; x>0,$$
 where $Y=-X$ and $\hat{U}$ is the
potential measure corresponding to $\hat{H}$.
 \end{lemma}
\begin{theorem}\label{ruin4-2}
\begin{itemize}
\item[(i)] Suppose $\pi_-$ is completely monotone
on $(0,\infty)$. Then, the probability of ruin $\psi$ is
completely monotone on $(0,\infty)$. In particular, $\psi\in C^{\infty}(0,\infty)$.
\item[(ii)] Suppose $\pi_-$ is  log-convex on $(0,\infty)$. Then,
\begin{itemize}
\item[(a)] $\psi$ is convex  on $(0,\infty)$;
\item[(b)] $\psi'$ is concave on $(0,\infty)$;
\item[(c)] if $X$ has no Gaussian component, then $\psi$ is twice continuously
differentiable except at finitely or countably many points  on
$(0,\infty)$, else $\psi \in C^2(0,\infty)$.
\end{itemize}
\end{itemize}
\end{theorem}
\noindent {\bf Proof.} \ We first prove (i). Since $\pi_-$ is completely monotone
on $(0,\infty)$, it follows from Lemma \ref{ruin4-1} that the tail $\Pi_{\cal
H}(x,\infty)$ of L\'evy measure $\Pi_{\cal H}$  is a complete
monotone function. Also, it follows from Lemma \ref{ladd-2} that  the
potential measure $U$ has a density $u$ which is completely monotone
on $(0,\infty)$. Thus,  the probability of ruin $\psi$ is completely
monotone on $(0,\infty)$ as $\psi(x)=\alpha U(x,\infty)$.

We now prove (ii). The log-convexity of $\pi_-$ implies the  log-convexity of
$\overline{\Pi}_Y^+$, and hence  $\overline{\Pi}_{\cal H}$ is
log-convex on $(0,\infty)$ due to Lemma \ref{ruin4-1} as
log-convexity is preserved under mixing. It follows from Lemma 3.1
that  the potential measure $U$ has a non-increasing and convex
density $u$. Thus, $\psi'=-\alpha u$ is non-decreasing and
concave on $(0,\infty)$, and hence (a) and (b) are proved.  Since a convex
function  on $(0,\infty)$ is differentiable except at finitely or
countably many points,  we see that $\psi$ is twice continuously differentiable except at finitely
or countably many points on $(0,\infty)$ if $X$ has no Gaussian component. On the other hand, if $X$ has a Gaussian
component, or equivalently, the drift of ascending ladder processes
$H$ strictly positive, then it follows from Lemma 3.1 that $u\in
C^1(0,\infty)$, and hence $\psi\in C^2(0,\infty)$. Therefore, (c) is proved.
\setcounter{equation}{0}
\section{ \normalsize Convexity of $h$}\label{scale}
\setcounter{equation}{0}
For  $h$ in (\ref{math-eq6}), define a barrier level
by
$$b^*=\sup\{b\ge 0: h'(b)\le h'(x) \;{\rm for\; all}\; x\ge 0\},$$
where $h'(0)$ is understood to be  the right-hand derivative at $0$.

For a spectrally negative L\'evy process, that is, in the case of
$\Pi\{(0,\infty)\}=0$, it was shown in Loeffen (2008) that the derivative of the $\delta$-scale function ${W^{(\delta)}}'(x)$ is
convex for $\delta > 0$ if $\Pi(x,\infty)$ is completely monotone.  This implies that there exists
an $a^*\ge 0$ such that $W^{(\delta)}$ is concave on $(0, a^*)$ and
convex on $(a^*,\infty)$.  Also, Kyprianou et al. (2010) showed
that if $\Pi(x,\infty)$  has a density on $(0,\infty)$ which is
non-increasing and log-convex, then for each $\delta \ge 0$, the scale
function $W^{(\delta)}(x)$ and its first derivative are convex beyond
some finite value of $x$.

Parallel to the results of Loeffen (2008) and Kyprianou et al. (2010) for spectrally negative L\'evy processes, we have the following  results.

\begin{theorem}\label{scale5-1}
\begin{itemize}
\item[(i)] Suppose $\pi_-$ is completely
monotone on $(0,\infty)$. Then, the derivative $h'(u)$
 is strictly convex on $(0,\infty)$ and
$h\in C^{\infty}(0,\infty)$;
\item[(ii)] Suppose $\pi_-$ is log-convex on $(0,\infty)$. Then, the $h$ and its derivative $h'$  are
strictly convex on $(b^*,\infty)$. Moreover, if $X$ has no Gaussian
component, $h$ is twice continuously differentiable except at
finitely or countably many points  on $(0,\infty)$, else $h\in C^2(0,\infty)$.
\end{itemize}
\end{theorem}
\noindent {\bf Proof.} \  Since $\pi_-$ is completely
monotone on $(0,\infty)$, we have $\tilde{\pi}_-$ is also completely
monotone on $(0,\infty)$, where $\tilde{\Pi}(dx)=\tilde{\pi}_-(-x)dx, x<0$. We can now apply  Theorem  \ref{ruin4-2} to deduce that
 the probability of ruin $\tilde{\psi}$ is completely monotone on $(0,\infty)$. In particular, $\tilde{\psi}\in C^{\infty}(0,\infty)$. It is easy to prove that $h'(u)$
 is strictly convex on $(0,\infty)$ and
$h\in C^{\infty}(0,\infty)$.  Hence, (i) is proved.

Let $\cal{\tilde H}$ ($\hat{\tilde{\cal H}}$) be the ascending (descending) ladder height process of
$\tilde{Y}=-\tilde{X}$. By Lemma  \ref{ruin4-1}, we have
$$\overline{\Pi}_{\cal{\tilde H}}(x)
=-\int_{-\infty}^0\hat{\tilde U}(dy)\overline{\Pi}_{\tilde
X}^{-}(x-y),\; \; x>0,$$
where  $\hat{\tilde U}$ is the renewal measure corresponding to $\hat{\tilde H}$. Then,
$$\Pi'_{\cal{\tilde H}}(x)
 =-e^{\rho(\delta)x}\int_{-\infty}^0 \hat{\tilde U}(dy)
 e^{-\rho(\delta)y}\pi_{-}(y-x)\equiv e^{\rho(\delta)x}\nu_+(x).$$
The assumption of log-convexity of $\pi_-$ implies that $\nu_+$ is
log-convex, and hence  $\Pi'_{\cal{\tilde H}}(x)$ is also
log-convex. It follows from Lemma 1 of Kyprianou and Rivero
 (2008) that the restriction of its potential measure to
$(0,\infty)$ of a subordinator with L\'evy density $\nu_+$ has a
non-increasing and convex density, say $f_{\delta}$. Also,  the
restriction of its potential measure to $(0,\infty)$ of a
subordinator with L\'evy density $\Pi'_{\cal{\tilde H}}(x)$ has a
non-increasing and convex density, say $h_{\delta}$. Moreover,
$h_{\delta}(x)=e^{\rho(\delta)x}f_{\delta}(x)$. Thus,
$\tilde{\psi}'(x)=-{\tilde \alpha} e^{\rho(\delta)x}f_{\delta}(x)$,
where $\tilde{\alpha}^{-1}=\int_0^{\infty}P(\tilde{H}_t<\infty)dt$.
Since $h(x)=[1-\tilde{\psi}(x)]e^{\rho(\delta) x}$, we have
$$h'(x)=\rho(\delta)h(x)+\alpha f_{\delta}(x),\; x>0.$$
This implies that $h'(x)$ tends to $\infty$ as $x$ tends to
$\infty$ as $\lim_{x\to\infty}h(x)=\infty$. Thus $b^*<\infty$.
Applying the same arguments as those in Kyprianou, Rivero and Song (2010), we can prove that $h$ and its derivative $h'$
are strictly convex on $(b^*,\infty)$.  Finally, the smoothness of $h$ is a direct consequence of Theorem  \ref{ruin4-2}. So, (ii) is proved.

\setcounter{equation}{0}
\section{\normalsize  Main results and proofs}\label{main}

We now present the main results of the paper about the
optimality of the barrier strategy $ \xi^{b^*}$ for the de
Finetti's dividend problem for general L\'evy  processes. This  is
a continuation of the work of Yuen and Yin (2011) in which a special L\'evy
process with both upward and downward jumps and a completely monotone density was considered.
\begin{theorem}\label{main-1} Suppose  that $\nu$ is a non-negative function on $(0,\infty)$
   which is sufficiently smooth and satisfies
  \begin{itemize}
  \item[(i)] $(\Gamma-\delta)\nu(x)\le 0,\ for\ almost\ every\ x>0;$
  \item[(ii)] $\nu \ is\ concave\ on \ (0,\infty);$
  \item[(iii)] $\nu'(x)\ge 1,\ x>0.$
   \end{itemize}
  Then, $\nu(x)\ge V_*(x)$.
 \end{theorem}
  \begin{theorem}\label{main-2}  Suppose that $V_b$ defined in (\ref{math-eq5}) is  sufficiently smooth
  and  satisfies
  \begin{itemize}
  \item[(i)] $V_b'(x)>1 \ for \ all\   x\in [0,b);$
  \item[(ii)] $(\Gamma-\delta)V_b(x)\le 0,\ for\ all\ x>b.$
  \end{itemize}
  Then, $V_b(x)= V_*(x)$. In particular, if $(\Gamma-\delta)V_{b^*}(x)\le 0$ for all $x>b^*$, then $V_{b^*}(x)= V_*(x)$.
\end{theorem}
 \begin{theorem}\label{main-3} Suppose that $\pi_-$ is completely
 monotone.  Then, $V_{b^*}(x)= V_*(x)$, that is, the barrier strategy at $b^*$ is the optimal strategy among all admissible strategies.
 \end{theorem}
\begin{theorem}\label{main-4} Suppose that $\pi_-$ is log-convex on $(0,\infty)$.
 Then, $V_{b^*}(x)= V_*(x)$, that is, the barrier strategy at $b^*$ is the optimal strategy among all admissible strategies.
 \end{theorem}
Before proving the main results, we give two lemmas which are similar to those in Loeffen (2008) for spectrally negative L\'evy processes.
 \begin{lemma}\label{main6-1}   Suppose that $h$  is  sufficiently smooth and convex in the interval $(b^*,\infty)$. Then, the following statements hold:
\begin{itemize}
\item[(i)] $b^*<\infty;$
\item[(ii)] $V_{b^*}'(x)\ge 1 \quad for\ x\in [0,b^*]$ and $V_{b^*}'(x)= V_{x}'(x)=1\ for\ x>b^*;$
\item[(iii)] $(\Gamma-\delta)V_{b^*}(x)=0 \quad for\ x\in (0,b^*).$
\end{itemize}
 \end{lemma}
\noindent{\bf Proof.} \ As $\lim_{x\to\infty}h'(x)=\infty$, we have (i).
For (ii), $V_{b^*}'(x)= h'(x)/h'(b^*)$ for $x\in [0,b^*]$; it follows from the definition of $b^*$ that $V_{b^*}(x)\ge 1$ for $ x\in [0,b^*]$; $V_{b^*}'(x)= V_{x}'(x)=1$ for $x>b^*$ because of $V_{b^*}(x)=x-b^*+V_{b^*}(b^*)$; and $V_x'(x)=1$ since $V_x(x)= h(x)/h'(x)$. Finally, (iii) is due to $(\Gamma-\delta)h(x)=0$ for $x\in (0,b^*)$ and (\ref{math-eq5}).
\begin{lemma}\label{main6-2} Suppose that $h$  is  sufficiently smooth and is convex in the interval $(b^*,\infty)$. Then, for $x>b^*$,
\begin{itemize}
\item[(i)]  $V_{b^*}''(x)=0\le V_x''(x-)$ \ if \ $\sigma\neq0$;
\item[(ii)] $V_{b^*}'(y)\ge V_x'(y),\  y\in [0,x]$;
\item[(iii)] $V_{b^*}(x)\ge V_x(x)$;
\item[(iv)] $(\Gamma-\delta)V_{b^*}(x)\le 0$.
\end{itemize}
\end{lemma}
\noindent{\bf Proof.} \ If  $\sigma\neq0$, $V_{b^*}''(x)=0$ is
clear. Also, since $h\in C^2(0,\infty)$ and is  convex in the interval
$(b^*,\infty)$, we have $V_x''(x-) = \lim_{y\uparrow x}V_x''(y)=\lim_{y\uparrow x} h''(y)/h'(x) \ge 0.$ Thus, (i) is proved.

For $y\in[0,b^*]$, by the definition of $b^*$, we have
 $$V_{b^*}'(y)-V_x'(y)=\frac{h'(y)}{h'(b^*)}-\frac{h'(y)}{h'(x)}\ge
 0.$$
On the other hand, for $y\in[b^*,x]$, by the convexity of $h$ on $(b^*,\infty)$, we
have
 $$V_{b^*}'(y)-V_x'(y)=1-\frac{h'(y)}{h'(x)}\ge
 0.$$
These give (ii).

Note that $V_{b^*}(b^*)= h(b^*)/h'(b^*) \ge h(b^*)/h'(x) = V_x(b^*)$ and that $(V_{b^*}-V_x)$ is
non-decreasing on $(b^*,\infty)$ because of (ii). Thus, $V_{b^*}(x)\ge V_x(x)$, that is, (iii) holds.

For $x>b^*$, $(\Gamma-\delta)V_{x}(x-) = \lim_{y\uparrow
x}(\Gamma-\delta)V_{x}(y)=0$. For $x \le b^*$, we have
\begin{eqnarray*}
(\Gamma-\delta)V_{b^*}(x)&=&(\Gamma-\delta)V_{b^*}(x)-(\Gamma-\delta)V_{x}(x-)\\
&=&\frac12\sigma^2(V_{b^*}''(x)-V_{x}''(x-))+a(V_{b^*}'(x)-V_{x}'(x))\\
&&+\int_{-\infty}^{\infty}(V_{b^*}(x+y)-V_{b^*}(x)-V_{b^*}'(x)y\text{\bf
1}_{\{|y|<1\}})\pi(y)dy\\
&&+\int_{-\infty}^{\infty}(V_{x}(x+y)-V_{x}(x)-V_{x}'(x)y\text{\bf
1}_{\{|y|<1\}})\pi(y)dy\\
&&-\delta(V_{b^*}(x)-V_{x}(x))\equiv I_1+I_2+I_3-I_4.
\end{eqnarray*}
Lemma  \ref{main6-1} (ii) and Lemma  \ref{main6-2} (i)  imply that
$I_1\le 0$, and Lemma \ref{main6-2} (iii) implies that $I_4\ge 0$. For $I_2+I_3$, we have
\begin{eqnarray*}
I_2+I_3&=&\int_{-\infty}^{\infty} ((V_{b^*}-V_x)(x+y)-(V_{b^*}-V_x)(x)-(V_{b^*}'-V_x')(x)y\text{\bf
1}_{\{|y|<1\}})\pi(y)dy\\
&=&\int_{-\infty}^{0} ((V_{b^*}-V_x)(x+y)-(V_{b^*}-V_x)(x)-(V_{b^*}'-V_x')(x)y\text{\bf
1}_{\{|y|<1\}})\pi(y)dy\\
&&+\int_{0}^{\infty} ((V_{b^*}-V_x)(x+y)-(V_{b^*}-V_x)(x)-(V_{b^*}'-V_x')(x)y\text{\bf
1}_{\{|y|<1\}})\pi(y)dy\\
&\equiv& J_1+J_2.
\end{eqnarray*}
Applying Lemma  \ref{main6-1} (ii) and Lemma  \ref{main6-2} (ii) yields $J_1\le 0$. For $y>0$, we obtain
$$(V_{b^*}-V_x)(x+y)=(V_{b^*}-V_x)(x)=x-b^*+\frac{h(b^*)}{h'(b^*)}-\frac{h(x)}{h'(x)},$$
which, together with Lemma  \ref{main6-2} (ii), imply that $J_2=0$. These prove (iv).

We now present the proofs of Theorems 6.1-6.4.

\noindent{\bf Proof of Theorem \ref{main-1}.} \ Define The jump measure of $X$ by
$$\mu^X=\mu^X(\omega,dt,dy)=\sum_s\text{\bf 1}_{\{\triangle X_s\neq
0\}}\delta_{(s,\triangle X_s)}(dt,dy),$$
and its compensator by $\upsilon=\upsilon(dt,dy)=dt\Pi(dy)$. Then, the L\'evy
decomposition (Protter, 1992, Theorem 42) gives
\begin{eqnarray*}
X_t&=&\sigma B_t+\int_{[0,t]\times {\Bbb R}}y\text{\bf
1}_{\{|y|<1\}}(\mu^X-\upsilon)+at+\int_{[0,t]\times {\Bbb
R}}(y-y\text{\bf 1}_{\{|y|<1\}})\mu^X\\
&\equiv& M_t+at+\sum_{0\le s\le t}\triangle X_s\text{\bf
1}_{\{|y|\ge 1\}},
\end{eqnarray*}
where $B=\{B_t\}_{t\ge 0}$ is a standard Brownian motion, and $M_t$ is a
martingale with $M_0=0$.

Note that $\nu$ is
smooth enough for an application of the appropriate version of
It\^{o}'s formula and the change of variables formula. In fact, if $X$
is of bounded variation, then $\nu\in C^1(0,\infty)$ and we are
allowed to use the change of variables formula (Theorem 31, Protter, 1992);
if $X$ has a Gaussian exponent, then $\nu\in C^2(0,\infty)$ and we
are allowed to use  It\^{o}'s formula (Theorem 32, Protter, 1992);
and if $X$ has unbounded variation and $\sigma=0$, then $\nu$ is twice
continuously differentiable almost everywhere but is not in
$C^2(0,\infty)$ and we can use the Meyer-It\^{o}'s formula (Theorem 70,
Protter, 1992) and product rule formula. In any cases,  for any appropriate localization
sequence of stopping times $\{t_n, n\ge 1\}$, we get under  $P_x$
\begin{eqnarray}
&&e^{-\delta  (t_n\wedge \tau^{\xi})}\nu({U}^{\xi}_{t_n\wedge \tau^{\xi}})-\nu({U}^{\xi}_0) \nonumber \\
&=&\int_{0}^{t_n\wedge \tau^{\xi}}
e^{-\delta s}dM_s^{\xi}+\int_0^{t_n\wedge \tau^{\xi}} e^{-\delta
s}(\Gamma-\delta)\nu({U}^{\xi}_{s-})ds \nonumber \\
&& \qquad +\sum_{s\le {t_n\wedge \tau^{\xi}}}\text{\bf 1}_{\{\triangle {L}_s^{\xi}>0\}}e^{-\delta s} \Bigg\{\nu({U}^{\xi}_{s-}+\triangle
X_s - \triangle {L}_s^{\xi})-\nu({U}^{\xi}_{s-} + \triangle X_s) \nonumber \\
&& \qquad \qquad + \nu'({U}^{\xi}_{s-}+\triangle X_s)\triangle
{L}_s^{\xi}\Bigg\} - \int_{0}^{t_n\wedge \tau^{\xi}} e^{-\delta
s}\nu'({U}^{\xi}_{s-})d{L}_s^{\xi}, \label{main-eq1}
\end{eqnarray}
where
\begin{eqnarray*}
M_t^{\xi}&=&\sum_{s\le t}\text{\bf 1}_{\{|\triangle
X_s|>0\}}\left\{\nu({U}^{\xi}_{s-}+\triangle
X_s)-\nu({U}^{\xi}_{s-})-\triangle X_s
\nu'({U}^{\xi}_{s-})\text{\bf 1}_{\{|\triangle X_s|\le 1\}}
\right\}\\
 &&-\int_0^t\int_{-\infty}^{\infty}\left\{\nu({U}^{\xi}_{s-}-y)
 -\nu({U}^{\xi}_{s-})+y\nu'({U}^{\xi}_{s-})\text{\bf 1}_{\{|y|\le
 1\}}\right\}\pi(y)dyds + \int_0^t\nu'({U}^{\xi}_{s-})dM_s
\end{eqnarray*}
is a local martingale.
The concavity of $\nu$ implies that $\nu(x)-\nu(y)+(x-y)\nu'(y)\le
0$ for any $x\le y$. Taking expectations on both sides of
(\ref{main-eq1}) and using the conditions (i)-(iii), we obtain
\begin{equation}
E_x (e^{-\delta (t_n\wedge \tau^{\xi})}\nu({U}^{\xi}_{t_n\wedge
\tau^{\xi}}))-\nu(x)\le -E_x\int_{0}^{t_n\wedge \tau^{\xi}} e^{-\delta s} d{L}_s^{\xi}. \label{main-eq2}
\end{equation}
Then, letting $n\to\infty$ in (\ref{main-eq2}) and recalling that $\xi$ is an arbitrary strategy in
$\Xi$, we get
$$\nu(x)\ge \sup_{ \xi\in\Xi}V_{\xi}(x)=V_*(x).$$
This ends the proof of Theorem \ref{main-1}.

\noindent{\bf Proof of Theorem  \ref{main-2}.} \
It follows from (\ref{math-eq5})  and the conditions (i) and (ii) that
 $(\Gamma-\delta)V_b(x)\le 0$ for $x\in (0,\infty)\setminus\{b\}$
and $V_b'(x)\ge 1$ for  $x>0$. Similar to (\ref{main-eq1}), one can show that
\begin{eqnarray}
&& e^{-\delta t}V_b({U}^{\xi}_t)-V_b({U}^{\xi}_0) \nonumber \\
&=& \int_0^t e^{-\delta t}dN_s^{\xi}+\int_0^t e^{-\delta
s}(\Gamma-\delta)V_b({U}^{\xi}_{s-})ds \nonumber \\
&& \qquad + \sum_{s\le t}\text{\bf 1}_{\{\triangle {L}_s^{\xi}>0\}}e^{-\delta
s}\left\{V_b({U}^{\xi}_{s-}+\triangle X_s-\triangle {L}_s^{\xi})-V_b({U}^{\xi}_{s-}+\triangle
X_s)\right\} \nonumber \\
&& \qquad - \int_{(0,t]} e^{-\delta s}V_b'({U}^{\xi}_{s-})d{L}_s^{\xi,c},\label{main-eq3}
\end{eqnarray}
where ${L}_s^{\xi,c}$ is the continuous part of
${L}_s^{\xi}$, and
\begin{eqnarray*}
N_t^{\xi}&=&\sum_{s\le t}\text{\bf 1}_{\{|\triangle
X_s|>0\}}\left\{V_b({U}^{\xi}_{s-}+\triangle
X_s)-V_b({U}^{\xi}_{s-})-\triangle X_s
V_b'({U}^{\xi}_{s-})\text{\bf 1}_{\{|\triangle X_s|\le 1\}}
\right\}\\
 &&-\int_0^t\int_{-\infty}^{\infty}\left\{V_b({U}^{\xi}_{s-}-y)
 -V_b({U}^{\xi}_{s-})+yV_b'({U}^{\xi}_{s-})\text{\bf 1}_{\{|y|\le
 1\}}\right\}\pi(y)dyds\\
&&+\int_0^tV_b'({U}^{\xi}_{s-})dM_s.
\end{eqnarray*}
Not that $P(\triangle{L}_s^{\xi}>0, \triangle X_s<0)=0$ and that ${U}^{\xi}_{s-}+\triangle X_s\ge b$ on $\{\triangle {L}_s^{\xi}>0, \triangle X_s>0\}$. Consequently, $V_b'({U}^{\xi}_{s-}+\triangle X_s)=1$, and hence
\begin{eqnarray*}
&\sum_{s\le t}\text{\bf 1}_{\{\triangle {L}_s^{\xi}>0\}}e^{-\delta
s}\left\{V_b({U}^{\xi}_{s-}+\triangle
X_s-\triangle {L}_s^{\xi})-V_b({U}^{\xi}_{s-}+\triangle
X_s)\right\}\\
&= -\sum_{s\le t}\text{\bf 1}_{\{\triangle {L}_s^{\xi}>0\}}e^{-\delta s}\triangle
{L}_s^{\xi}.
\end{eqnarray*}
Also, for any appropriate localization sequence of stopping times $\{t_n, n\ge 1\}$, we have
\begin{equation}
E_x (e^{-\delta (t_n\wedge
\tau^{\xi})}V_b({U}^{\xi}_{t_n\wedge
\tau^{\xi}}))-E_x V_b({U}^{\xi}_0)\le -E_x\int_{[0, t_n\wedge
\tau^{\xi}]} e^{-\delta s} d {L}_s^{\xi}.\label{main-eq4}
\end{equation}
Letting $n\to\infty$ in (\ref{main-eq4}) yields
$$V_b(x)\ge \sup_{ \xi\in\Xi}V_{\xi}(x)=V_*(x).$$
However, $$V_b(x)\le \sup_{ \xi\in\Xi}V_{\xi}(x)=V_*(x).$$
This ends the proof of Theorem  \ref{main-2}.

\noindent{\bf Proof of Theorem  \ref{main-3}.} \ If $\pi_-$ is completely
monotone, it follows from  Theorem  \ref{scale5-1} (i) that $h'(x)$
is strictly convex on $(0,\infty)$. Then, $V_{b^*}$ is concave on $(0,\infty)$ because of (\ref{math-eq5}). From Lemma \ref{main6-1}
(ii) and (iii) and  Lemma  \ref{main6-2}  (iv), we see that the
conditions in Theorem  \ref{main-1} are satisfied. Thus, $V_b(x)\ge
V_*(x)$. Consequently, $V_b(x)= V_*(x)$ and the proof is complete.

\noindent{\bf Proof of Theorem   \ref{main-4}.} \ If $\pi_-$ is  log-convex
on $(0,\infty)$, it follows from  Theorem  \ref{scale5-1} (ii) that
$h(x)$ is strictly convex on $(b^*,\infty)$. Then, applying Lemma  \ref{main6-2} (iv) gives $(\Gamma-\delta)V_{b^*}(x)\le 0$
for all $x>b^*$. The result follows from Theorem \ref{main-2}.

\vskip0.3cm

\noindent{\bf Acknowledgements}\; The authors are grateful to Dr. Yang Xuewei of Nankai University for useful discussions.
The research of Chuancun Yin was supported by
the National Natural Science Foundation of China (No. 11171179) and
the Research Fund for the Doctoral Program of Higher Education of
China (No. 20133705110002). The research of Kam C. Yuen was supported by a
university research grant of the University of Hong Kong.

\end{document}